# THE CRITICAL ISING MODEL ON TREES, CONCAVE RECURSIONS AND NONLINEAR CAPACITY

By Robin Pemantle[1] and Yuval Peres[2]

*University of Pennsylvania and Microsoft Research*

We consider the Ising model on a general tree under various boundary conditions: all plus, free and spin-glass. In each case, we determine when the root is influenced by the boundary values in the limit as the boundary recedes to infinity. We obtain exact capacity criteria that govern behavior at critical temperatures. For plus boundary conditions, an $L^3$ capacity arises. In particular, on a spherically symmetric tree that has $n^\alpha b^n$ vertices at level $n$ (up to bounded factors), we prove that there is a unique Gibbs measure for the ferromagnetic Ising model at the relevant critical temperature if and only if $\alpha \le 1/2$. Our proofs are based on a new link between nonlinear recursions on trees and $L^p$ capacities.

**1. Introduction.** Let $T$ be a finite rooted tree of depth $N$. Let $|v|$ denote the distance from a vertex $v \in V(T)$ to the root $o$, and write $v \to w$ if $v$ is the parent of $w$, that is, the neighbor of $w$ closer to the root than $w$. Consider the space $\Omega = \Omega(T) = \{+1, -1\}^{V(T)}$ of configurations on the vertices of $T$. For each $w \ne o$ there is a unique edge $vw$ with $v \to w$; let $J_w = J(vw)$ be a positive number so that $\{J_w : o \ne w \in V(T)\}$ is a fixed set of interaction strengths on the edges of $T$. We assume throughout that the interaction strengths are bounded;

$$0 < J_{\min} \le J_v \le J_{\max} \qquad \forall v \in V(T), v \ne o. \tag{1.1}$$

This assumption loses little generality (see the end of Section 4). Fix an inverse temperature $\beta$ and define the weight of a configuration $\eta \in \Omega$ to be the following product over all pairs of neighboring vertices:

$$W(\eta) = \prod_{v \to w} \exp(\beta J_w \eta(v)\eta(w)).$$

Received May 2006; revised April 2009.
[1]Supported in part by NSF Grant DMS-01-03635.
[2]Supported in part by NSF Grants DMS-01-04073 and DMS-02-44479.

The Ising model under various boundary conditions can be obtained by restricting to suitable subsets of $\Omega$ and assigning probabilities proportional to $W$. Our aim in this paper is to pinpoint the locations of the phase transitions that occur in these models as $N \to \infty$. In each case the critical temperature for phase transitions to occur is known. We refine these results by giving sharp criteria for the existence of a phase transition in terms of capacities.

**2. Main results.** Let $T$ be any tree rooted at a vertex $o$, and let $\partial T$ denote the set of maximal paths oriented away from the root; these are either infinite or end at a leaf of $T$. For finite trees, we identify $\partial T$ with the set of leaves in $T$ different from $o$. For infinite trees, we assume there are no leaves (except possibly $o$), so all paths in $\partial T$ are infinite. Let $\{R(e) : e \in E(T)\}$ be a set of *resistances* (nonnegative numbers) assigned to the edges of $T$. Let $\mu$ be a flow on $T$, that is, a nonnegative function on $E(T)$ such that at every vertex (except for the root and the leaves) inflow equals outflow; whenever $v \to w$ ($v$ is a parent of $w$) and $w$ is not a leaf, we have $\mu(vw) = \sum_{y\,:\,w\to y} \mu(wy)$. Such a flow $\mu$ can be identified with a positive finite measure on $\partial T$ where $\mu(e)$ is the measure of the set of paths in $\partial T$ that traverse $e$. The total mass of this measure is the outflow from the root, $|\mu| := \sum_{y\,:\,o\to y} \mu(oy)$. Fix $p > 1$ and set $s = p - 1$. For $y \in \partial T$, define

$$V_\mu(y) := \sum_{e\in y} (\mu(e)R(e))^s; \tag{2.1}$$

$$V(\mu) := \sup\{V_\mu(y) : y \in \partial T\}; \tag{2.2}$$

$$\operatorname{cap}_p(T) := \sup\{|\mu| : \mu \text{ a flow on } T \text{ with } V(\mu) = 1\}. \tag{2.3}$$

These capacities have been studied on more general networks as part of discrete nonlinear potential theory [see, e.g., Murakami and Yamasaki (1992), Soardi (1993, 1994) and the references therein]. However, all the properties of $\operatorname{cap}_p$ that we will use follow readily from the definition. We note that $\operatorname{cap}_2(T)$ reduces to the electrical conductance between $o$ and $\partial T$. We also observe that if the tree $T$ and the resistances are spherically symmetric [i.e., the degree of every vertex $w$ depends only on $|w|$, and similarly for the resistance $R(vw)$ between $w$ and its parent $v$], then among all flows $\mu$ with the same mass $|\mu|$, the equally splitting flow minimizes $V(\mu)$. To see this, given any other $\mu$, choose a path from $o$ to $\partial T$ by maximizing $\mu$ at every step.

For $T$ finite, let $\mathbf{P}$ denote the probability measure on $\Omega(T)$ proportional to $W$;

$$\mathbf{P}(\eta) = \frac{W(\eta)}{\sum_{\xi\in\Omega} W(\xi)}.$$

This is a ferromagnetic Ising model with no external field and free boundary conditions. There is another construction of the measure $\mathbf{P}$ as a tree-indexed

Markov chain. To the edge leading to a vertex $v$ from its parent, assign the positive *bias*,

$$\theta_v = \frac{e^{\beta J_v} - e^{-\beta J_v}}{e^{\beta J_v} + e^{-\beta J_v}} = \tanh(\beta J_v). \tag{2.4}$$

Let the spin $\eta(o)$ at the root take the values $\pm 1$ with probability $1/2$ each. Conditional on $\eta(o)$, let the spin at every other vertex $v$ be determined recursively, by copying the sign at the parent with probability $(1+\theta_v)/2$ and reversing sign with probability $(1-\theta_v)/2$. When $J_v = J$ does not depend on $v$, we write $\theta$ for the common bias.

Now suppose $T$ is an infinite, locally finite tree, rooted at o, and let $T^{(N)}$ be the induced finite subgraph of $T$ with vertices $\{v \in V(T) : |v| \le N\}$. Letting $\mathbf{P}^{(N)}$ be the free-boundary Ising measure on $T^{(N)}$, we ask about $\mathbf{P}^{(N)}(\eta(0) = +1 | \eta(v) : v \in \partial T^{(N)})$. In particular, this converges in probability to $1/2$ if and only if the free boundary Gibbs measure on $T$ is extremal [see Georgii (1988)]. The question of extremality of the Gibbs measure with free boundary on regular trees was settled by Bleher, Ruiz and Zagrebnov (1995) [see also Ioffe (1996a) for an elegant alternative proof].

The same question for general trees was solved by Ioffe (1996b) and Evans et al. (2000) where the critical value is computed for an arbitrary tree. However, the question of extremality at the critical temperature was left open. In this paper we settle the critical case by showing that zero $L^2$ capacity (with respect to certain resistances) implies extremality. For vertices $y, w$ of $T$, write $y \le w$ if $y$ is on the path from the root $o$ to $w$. If $y \le w$ and $y \ne w$ write $y < w$. In particular, $o < y$ for every vertex $y \ne o$. We prove:

THEOREM 2.1. *Let $T$ be an infinite, locally finite tree, rooted at $o$, with no leaves except possibly at $o$ and interaction strengths $J_v$ satisfying (1.1). For vertices $y, w$, write $y \le w$ if $y$ is on the path from $o$ to $w$. Assign to each edge $e = vw$ with $v \to w$, the resistance*

$$R_w := R(e) := \prod_{o < y \le w} (\tanh \beta J_y)^{-2}. \tag{2.5}$$

*Then the free boundary Gibbs measure at inverse temperature $\beta$ is extremal if and only if* $\mathrm{cap}_2(T) = 0$.

One direction of this theorem (that extremality implies zero capacity) was already proved in Evans et al. (2000).

*Plus boundary conditions.* Consider $T$ finite again. Let $\Omega_+ = \Omega_+(T) \subset \Omega(T)$ be the set of configurations with $\eta(v) = +1$ for $v \in \partial T$. Then the probability measure $\mathbf{P}^+$ on $\Omega_+$ defined by

$$\mathbf{P}^+(\eta) = \frac{W(\eta)}{\sum_{\xi \in \Omega_+} W(\xi)}$$

is the Ising model with plus boundary conditions and no external field.

The critical value of the interaction strength here has long been known for regular trees [see Preston (1974, 1976)]. Lyons (1989) computes the critical temperature for general trees and allows the interaction strengths to vary as well. We refine the known results by determining what happens at criticality. The sharp criterion turns out to involve an "$L^3$-capacity." We prove:

THEOREM 2.2. *Let $T$ be any infinite, locally finite tree rooted at $o$ and having no leaves except possibly $o$. Let $\{J_v\}$ be bounded interaction strengths, that is, satisfying (1.1), and assign resistance $R_v = \prod_{o<y\le v}(\tanh(\beta J_y))^{-1}$ to the edge between $v$ and its parent. Then the decreasing limit*

$$\lim_{N\to\infty} \mathbf{P}^{(N,+)}(\eta(o)=+1)$$

*is equal to $1/2$ if and only if $\operatorname{cap}_3(T)=0$.*

Here $\mathbf{P}^{(N,+)}$ is the measure on configurations on the first $N$ levels of $T$ with plus boundary conditions imposed at level $N$.

For ease of reading, we state the result more explicitly in the special case of spherically symmetric trees and when the interaction strength is constant.

COROLLARY 2.3. *Under the hypotheses of Theorem 2.2, assume spherical symmetry as well; $\theta_v = \theta_{|v|}$ and $\deg(v) = d_{|v|}$ depend only on $|v|$. Then there are multiple Gibbs states if and only if*

$$\sum_{n\ge 1}\prod_{i=1}^{n}(d_i\theta_i)^{-2} < \infty. \tag{2.6}$$

*In particular, for a spherically symmetric tree $T$, suppose that the level cardinalities satisfy*

$$|T_n| \asymp \theta^{-n} n^{\alpha}. \tag{2.7}$$

*Then there is a unique Gibbs state for the Ising model at criticality if and only if $\alpha \le 1/2$.*

Note that for $T$ satisfying (2.7), endowed with edge resistances $\theta^{-n}$ at level $n$, the standard $L^2$ capacity of $T$ is zero as long as $\alpha \le 1$.

COROLLARY 2.4. *Suppose that $J_v \equiv J$ is constant, and let $\theta := \tanh(\beta J)$. Then phase transition occurs with plus boundary conditions if and only if $\operatorname{cap}_3(T) > 0$ with resistances $\theta^{-n}$ at distance $n$ from the root. If $T$ is spherically symmetric, this is equivalent to*

$$\sum_{n\ge 1}\theta^{-2n}|T_n|^{-2} < \infty.$$

*(The last statement is also a special case of the previous corollary.)*

*Spin-glass boundary conditions.* For a tree $T$ of depth $N$, define a measure $\mathbf{P}^{\mathrm{sg}}$ on $\Omega(T)$ by making the signs $\eta(v)$ for $v \in \partial T$ i.i.d. fair coin flips and requiring that the measure be proportional to $W$ conditionally on the values on $\partial T$:

$$\mathbf{P}^{\mathrm{sg}}(\eta) = 2^{-|\partial T|} \frac{W(\eta)}{\sum_{\xi|_{\partial T} = \eta|_{\partial T}} W(\xi)}.$$

This is equivalent to the following spin-glass model considered by Chayes et al. (1986): the Hamiltonian has interactions of a fixed magnitude, and no external field; the signs of the interactions are determined by i.i.d. fair coin flips, and the boundary conditions are fixed and known (e.g., they are all plus). The question is whether, conditional upon the signs of the interactions, the sign at the root is influenced at all by the boundary values in the limit as $N \to \infty$. A critical interaction strength is given in Chayes et al. (1986) for regular trees; we improve this to the case of general trees and settle what happens at the critical case. The result is a standard (i.e., $L^2$) capacity criterion, exactly equal to the criterion for the case of a free boundary.

THEOREM 2.5. *Let $T$ be an infinite, locally finite tree, rooted at $o$, with no leaves (except possibly $o$) and interaction strengths $J_v$ satisfying (1.1). Assign resistances*

$$R_v = \prod_{0 < y \le v} (\tanh(\beta J_y))^{-2}.$$

*Then $\mathbf{P}^{\mathrm{sg}}(\eta(o) = 1 | \eta|_{\partial T^{(N)}}) \to 1/2$ in law under the spin-glass measure if and only if $\mathrm{cap}_2(T) = 0$.*

*Recursions for the log-likelihood.* Let $x_v$ denote the log-likelihood ratio of having spin 1 versus $-1$ at $v$, given the boundary. The method in the plus boundary case is to show that $\{x_v : v \in V(T)\}$ satisfy a recursion of the form

$$x_v = \sum_{v \to w} f_w(x_w). \tag{2.8}$$

This reduces the problem to the question of whether, on a given infinite tree, this recursion has a nonzero solution. We give a general solution to this problem, recursively establishing a set of inequalities relating solutions and sub-solutions of these equations to generalized capacities. In the cases of free and spin-glass boundary conditions, the log likelihood ratios are random variables $\{X_v : v \in V(T)\}$, and we obtain versions of (2.8) for certain moments $\{m_v\}$ of $\{X_v\}$.

The rest of the paper is organized as follows. The next section focuses entirely on the deterministic aspect of the problem, namely, when the recursion (2.8) has a nontrivial solution or sub-solution. The theorems in this section are broad enough to handle the recursions arising from the three types of boundary conditions in the Ising model. Then we spend one section on each of the three models and conclude with some questions.

**3. Recursions on trees and potential theory.** Let $T$ be any locally finite rooted tree and let $\{f_v : v \in V(T)\}$ be a collection of nonnegative functions indexed by the vertices of $T$. We are interested in whether the simultaneous inequalities,

$$x_v \le \sum_{v \to w} f_w(x_w), \tag{3.1}$$

have any nonzero solutions. A special case of interest is when $f_v \equiv f$ does not depend on $v$. Our characterization is in terms of generalized capacities, which we defined in (2.3).

Fix $p > 1$ and let $s = p - 1$. We quote several easy and well-known consequences of the definition of capacity:

(i) The supremum in the definition (2.3) of $\mathrm{cap}_p$ is achieved if the set of measures of bounded potential is nonempty. [Clear by lower semi-continuity of $V(\mu)$.]

(ii) Joining several trees at the root sums their capacities.

(iii) Multiplying all resistances by $\alpha$ decreases capacity by a factor of $\alpha$.

(iv) A single edge of resistance $R$ connected in series to the root of a tree $T$ yields a tree of capacity,

$$\frac{\mathrm{cap}_p(T)}{(1 + R^s \, \mathrm{cap}_p(T)^s)^{1/s}}.$$

To see (iv), observe that there is a one-to-one correspondence between flows $\mu$ from the root to the boundary in $T$ and flows $\mu_R$ in the enhanced tree, such that $|\mu_R| = |\mu|$ and $V(\mu_R) = V(\mu) + R^s |\mu|^s$.

These facts yield the following lemma, which we will need below. Denote by $T(v)$ the subtree of $T$ consisting of $v$ and all vertices that are separated from $o$ by $v$.

Lemma 3.1. *Fix $p > 1$ and $s = p - 1$. For any vertex $v$, define*

$$\phi(v) := R_v \, \mathrm{cap}_p(T(v)),$$

*where $R_o = 1$ by convention. [In particular, $\phi(v) = R_v$ if $v$ is a leaf.] Then for any vertex $v$,*

$$\phi(v) = \sum_{v \to w} \frac{(R_v/R_w)\phi(w)}{(1 + \phi(w)^s)^{1/s}}.$$

PROOF. If $w \neq o$, let $T'(w)$ be the tree rooted at the parent of $w$ consisting of $T(w)$ plus the edge between $w$ and its parent. Then

$$\begin{aligned}\phi(v) = R_v \operatorname{cap}_p(T(v)) &= \sum_{v \to w} R_v \operatorname{cap}_p(T'(w)) \\ &= \sum_{v \to w} (R_v/R_w) \frac{R_w \operatorname{cap}_p(T(w))}{(1 + R(w)^s \operatorname{cap}_p(T(w))^s)^{1/s}},\end{aligned}$$

which gives the desired expression. □

We now relate these computations to the system (3.1). In the following theorem, $f(\infty)$ denotes $\liminf_{x \to \infty} f(x)$ and $s$ denotes $p-1$.

THEOREM 3.2. *Let $T$ be finite. Suppose that there exist $\kappa_1 > 0$, $p = 1 + s > 1$ and a collection of positive constants $\{a_v : v \in V(T)\}$ such that for every $v \in V(T)$ and $x \geq 0$,*

$$f_v(x) \leq \frac{a_v x}{(1 + (\kappa_1 x)^s)^{1/s}}. \tag{3.2}$$

*Then any solution to the system*

$$x_v = \sum_{v \to w} f_w(x_w) \qquad \text{with } x_w = \infty \text{ when } w \text{ is a leaf,} \tag{3.3}$$

*satisfies*

$$x_o \leq \frac{\operatorname{cap}_p(T)}{\kappa_1}, \tag{3.4}$$

*where the resistances are given by*

$$R_v = \prod_{0 \leq y \leq v} a_y^{-1}. \tag{3.5}$$

*Similarly, if (3.3) holds and*

$$\frac{a_v x}{(1 + (\kappa_2 x)^s)^{1/s}} \leq f_v(x), \tag{3.6}$$

*then*

$$\frac{\operatorname{cap}_p(T)}{\kappa_2} \leq x_o. \tag{3.7}$$

PROOF. We first prove that (3.3) and (3.6) imply (3.7). Let $g(v) = R_v \times \operatorname{cap}_p(T(v))/\kappa_2$, with $g(v) = \infty$ if $v \neq o$ is a leaf. We show by induction that $g(v) \leq x_v$ for all $v$. If $v$ is a leaf, this is true by definition. Assume $v$ is not a leaf and, by induction, that $g(w) \le x_w$ for all $v \to w$. Applying the previous lemma gives

$$g(v) = \sum_{v \to w} \frac{(R_v/R_w)g(w)}{(1 + (\kappa_2 g(w))^s)^{1/s}}.$$

Note that $R_v/R_w = a_w$ when $v \to w$. By monotonicity of $x \mapsto x/(1+(cx)^s)^{1/s}$, and the induction hypothesis,

$$g(v) \le \sum_{v \to w} \frac{a_w x_w}{(1 + (\kappa_2 x_w)^s)^{1/s}}.$$

This is, at most, $\sum_{v \to w} f_w(x_w)$ by the assumption (3.6), finishing the induction.

If we assume (3.2) instead of (3.6), an analogous induction yields $x_v \le G(v)$ for all $v$ where $G(v) = R_v \operatorname{cap}_p(T(v))/\kappa_1$. Setting $v = o$ now recovers the statement of the theorem. □

With regard to sub-solutions, that is, to the system of inequalities (3.1), we have the following immediate corollary, used in Section 4 to analyze Ising models with free boundaries.

Corollary 3.3. *Under the hypothesis (3.2), any solution to*

$$x_v \le \sum_{v \to w} f_w(x_w)$$

*satisfies* $x_o \le \operatorname{cap}_p(T)/\kappa_1$.

Although the estimate for finite trees in Theorem 3.2 is the most useful, the following corollary for infinite trees is more elegant. The corollary follows directly from the fact that $\operatorname{cap}_p(T)$ is the decreasing limit of $\operatorname{cap}_p(T^{(N)})$, so we omit the details.

Corollary 3.4. (a) *Let $T$ be infinite and locally finite, having no leaves except possibly the root. Assign resistances according to (3.5). If $f$ satisfies (3.2) for all $v \in V(T), x \ge 0$, then any solution $\{x_v\}$ of (3.1) satisfies $x_o \le \operatorname{cap}_p(T)/\kappa_1$. In particular, if $\operatorname{cap}_p(T) = 0$ and (3.2) holds, then there are no nontrivial solutions to (3.1) on $T$.*

(b) *Conversely, if $\operatorname{cap}_p(T) > 0$ and $f$ satisfies (3.6) for all $v \in V(T)$ and $x \ge 0$, then there is a solution of (3.1) with the property that $x_o \ge \operatorname{cap}_p(T)/\kappa_2$. This solution is given by $x_v = R_v \operatorname{cap}_p(T(v))/\kappa_2$ for all $v$.*

To see the value in what we have proved, we turn to some special cases. Recall that we denote $f(\infty) = \liminf_{x \to \infty} f(x)$.

COROLLARY 3.5. *Suppose that an increasing bounded function* $f:[0,\infty)\to[0,\infty)$ *satisfies:*

(i) $f(x)=ax-\Theta(x^p)$ *near 0 for some* $p>1$;
(ii) $0<f(x)<ax$ *for all* $x>0$.

*Then there is a nontrivial sub-solution* $x_v\le\sum_{x\to w}f(x_w)$ *on the vertices of* $T$ *if and only if* $\mathrm{cap}_p(T)>0$ *with resistances* $a^{-n}$ *at distance* $n$ *from the root.*

REMARKS.

- Assumption (ii) above follows from (i) if $f$ is concave and $f(x)>0$ for all $x>0$.
- Denote by $|e|$ the level of an edge $e$, so edges adjacent to $o$ have $|e|=1$. The *branching number* $\mathrm{br}(T)$ of an infinite tree $T$ was defined by Lyons (1990, 1992) as the infimum of the $\lambda$ such that $T$ admits a nonzero flow $\mu$ that satisfies $\mu(e)\le\lambda^{-|e|}$ for all edges $e$ of $T$. Suppose we assign resistance $R(e)=a^{-|e|}$ to every edge $e$ of $T$. If an infinite tree $T$ has $\mathrm{br}(T)<a^{-1}$ then any positive flow $\mu$ on $T$ must satisfy $\mu(e)\ge(a+\delta)^{|e|}$ for some $\delta>0$ and infinitely many edges $e$, whence $V(\mu)=\infty$. Thus
$$\mathrm{br}(T)<a^{-1}\quad\Longrightarrow\quad\mathrm{cap}_p(T)=0\qquad\forall p>0.$$
Conversely,
$$\mathrm{br}(T)>a^{-1}\quad\Longrightarrow\quad\mathrm{cap}_p(T)>0\qquad\forall p>0,$$
since under this assumption, $T$ admits a flow $\mu$ with $\mu(e)\le(a-\delta)^{|e|}$ for some $\delta>0$ and all edges $e$.
- Lyons (1990, 1992) proved that for Bernoulli percolation on a tree $T$ with retention probability $a$ for each edge, the probability that the root is in an infinite cluster satisfies $\mathbf{P}[o\longleftrightarrow\partial T]>0$ iff $\mathrm{cap}_2(T)>0$ where the resistance of an edge $e$ is $a^{-|e|}$. One of the proofs Lyons gave was recursive, and it was refined by Marchal (1998). This result is covered by our framework (though not with the optimal constants);

  Define $x_v:=-\log(1-\mathbf{P}[v\longleftrightarrow\partial T(v)])$ and rewrite the identity
$$1-\mathbf{P}[v\longleftrightarrow\partial T(v)]=\prod_{\{w\,:\,v\to w\}}(1-a\mathbf{P}[w\longleftrightarrow\partial T(w)])$$
in the form
$$e^{-x_v}=\prod_{w\,:\,v\to w}(1-a(1-e^{-x_w})),$$
that is, $x_v=\sum_{v\to w}f(x_w)$ where
$$f(x)=-\log[1-a(1-e^{-x})].$$

It is easy to check that $f(x) = ax - \Theta(x^2)$ near 0 and $f$ is concave, so it satisfies the hypotheses of Corollary 3.5, whence the claimed equivalence for percolation follows.

COROLLARY 3.6. *Suppose that $T$, an infinite, locally finite, leafless tree, is spherically symmetric, meaning that the degree of $v$ depends only on $v$. Suppose $f_v = f_{|v|}$ depends only on $|v|$ as well. Assume the inequalities (3.2) and (3.6). Then there is a nonzero solution to $x_v \le \sum_{v\to w} f_w(x_w)$ if and only if*

$$\sum_{n=1}^{\infty}\prod_{j=1}^{n}\frac{1}{(a_j d_j)^s} < \infty.$$

We use this in the next section with $s = 2$ to obtain an exact summability criterion for phase transition of the Ising model with plus boundary conditions on an arbitrary spherically-symmetric tree. This refines the work of Lyons (1989), who computed the critical value in terms of the branching number but did not settle the behavior at criticality.

**4. Plus boundary conditions.** In this section $T$ is an infinite tree with no leaves except possibly the root, and $T^{(N)}$ denotes the truncation to distance at most $N$ from the root. We fix interaction strengths $\{J_w : o \ne w \in V(T)\}$ satisfying (1.1), set $\theta_v = \tanh(\beta J_v)$ and consider the family of measures $\mathbf{P}^{(N,+)}$ on the space $\Omega_+(T^{(N)})$ of $\pm 1$ configurations on $T^{(N)}$ with plus boundary conditions. The goal is to determine whether $\mathbf{P}^{(N,+)}(\eta(0) = +1)$ converges to $1/2$ or is bounded below by $1/2 + \varepsilon$ as $N \to \infty$. This is accomplished in the following theorem, already stated in the Introduction.

THEOREM 4.2. *Let $T$ be any infinite, locally finite tree rooted at $o$ and having no leaves except possibly $o$. Let $\{J_v\}$ be bounded interaction strengths, that is, satisfying (1.1), and assign resistances $R_v = \prod_{0<y\le v}(\tanh(\beta J_y))^{-1}$ as in (2.5). Then the decreasing limit,*

$$\lim_{N\to\infty}\mathbf{P}^{(N,+)}(\eta(0) = +1),$$

*is equal to $1/2$ if and only if $\mathrm{cap}_3(T) = 0$.*

The key to the proof of Theorem 2.2 and to the main results in each of the next two sections is the following recursive likelihood computation. For any tree denote by $T(v)$ the subtree rooted at $v$ so that, for $|v| \le N$, the tree $T^{(N)}(v)$ has vertex set $\{w \in V(T) : v \le w, |w| \le N\}$. Consider a boundary configuration $\xi : \partial T^{(N)} \to \{\pm 1\}$ and let $\mathbf{P}^{\xi}$ denote the Ising measure with boundary condition $\xi$. Furthermore, let $\mathbf{P}_v^{(N,\xi)}$ denote the Ising measure on $T^{(N)}(v)$ whose boundary condition is $\xi_{|\partial T^{(N)}(v)}$.

LEMMA 4.1. *For each* $v \neq o$ *let* $\theta_v = \tanh(\beta J_v) \in [0,1)$. *Let*

$$x_v^{(N)} = x_v^{(N,\xi)} = \log\left[\frac{\mathbf{P}_v^{(N,\xi)}(\eta(v)=+1)}{\mathbf{P}_v^{(N,\xi)}(\eta(v)=-1)}\right]$$

*be the log-likelihood ratio at the root given the boundary. Then for* $|v| < N$,

$$x_v^{(N)} = \sum_{v\to w} f_w(x_w^{(N)}),$$

*where for* $\theta \in [0,1)$ *and* $w \in V(T)$, *we denote*

$$f_\theta(x) := \log\left[\frac{\cosh(x/2)+\theta\sinh(x/2)}{\cosh(x/2)-\theta\sinh(x/2)}\right] \quad \textit{and} \quad f_w(x) := f_{\theta_w}(x). \tag{4.1}$$

This lemma is well known; we include its proof for the convenience of the reader.

PROOF OF LEMMA 4.1. Let $\eta$ be a configuration on $T^{(N)}(v)$. If $|v| < N$ then for each child $w$ of $v$, let $\eta_w$ be the restriction of $\eta$ to the subtree $T^{(N)}(w)$. We may then write

$$W(\eta) = \prod_{v\to w} W(\eta_w)\exp(\eta(v)\eta(w)\beta J_w).$$

Writing $Z_v$ for the normalizing factor, we have

$$\mathbf{P}_v^{(N,\xi)}(\eta(v)=+1) = Z_v^{-1}\prod_{v\to w}\sum_{\eta_w : T^{(N)}(w)\to\{\pm 1\}} W(\eta_w)\exp(\eta(w)\beta J_w),$$

which equals

$$Z_v^{-1}\prod_{v\to w}[e^{\beta J_w} Z_w \mathbf{P}_w^{(N,\xi)}(\eta(w)=1) + e^{-\beta J_w} Z_w \mathbf{P}_w^{(N,\xi)}(\eta(w)=-1)]. \tag{4.2}$$

Similarly, $\mathbf{P}_v^{(N,\xi)}(\eta(v)=-1)$ equals

$$Z_v^{-1}\prod_{v\to w}[e^{-\beta J_w} Z_w \mathbf{P}_w^{(N,\xi)}(\eta(w)=+1) + e^{\beta J_w} Z_w \mathbf{P}_w^{(N,\xi)}(\eta(w)=-1)]. \tag{4.3}$$

Divide (4.2) and (4.3) by $\prod_{v\to w} Z_w \mathbf{P}_w^{(N,\xi)}(\eta(w)=-1)$ and then consider their ratio;

$$\begin{aligned}\frac{\mathbf{P}_v^{(N,\xi)}(\eta(v)=+1)}{\mathbf{P}_v^{(N,\xi)}(\eta(v)=-1)} &= \prod_{v\to w}\frac{e^{(\beta J_w + x_w^{(N)})} + e^{-\beta J_w}}{e^{(-\beta J_w + x_w^{(N)})} + e^{\beta J_w}} \\ &= \prod_{v\to w}\frac{\cosh(\beta J_w)(e^{x_w^{(N)}}+1) + \sinh(\beta J_w)(e^{x_w^{(N)}}-1)}{\cosh(\beta J_w)(e^{x_w^{(N)}}+1) - \sinh(\beta J_w)(e^{x_w^{(N)}}-1)}.\end{aligned}$$

Next, divide numerator and denominator by $\cosh(\beta J_w)$ and recall that $\tanh(\beta J_w) = \theta_w$. It follows that the log of the likelihood ratio above satisfy

$$x_v^{(N)} = \sum_{v \to w} \log \frac{e^{x_w^{(N)}} + 1 + \theta_w (e^{x_w^{(N)}} - 1)}{e^{x_w^{(N)}} + 1 - \theta_w (e^{x_w^{(N)}} - 1)}. \tag{4.4}$$

Finally, divide numerator and denominator by $e^{x_w^{(N)}/2}$ to complete the proof. $\square$

We will need some basic properties of the functions $f_\theta$ defined in (4.1).

LEMMA 4.2. *For $\theta > 0$, the function*

$$f_\theta(x) := \log \left[ \frac{\cosh(x/2) + \theta \sinh(x/2)}{\cosh(x/2) - \theta \sinh(x/2)} \right]$$

*is an increasing odd function of $x \in \mathbf{R}$ which is concave for $x > 0$. Moreover, for any compact interval $I \subset (0, \infty)$, the inequality*

$$\frac{\theta x}{(1 + \kappa_2 x^2)^{1/2}} \le f_\theta(x) \le \frac{\theta x}{(1 + \kappa_1 x^2)^{1/2}} \tag{4.5}$$

*holds for all $x > 0$ and $\theta \in I$ where the constants $\kappa_2 \ge \kappa_1 > 0$ depend only on $I$.*

PROOF. First, we differentiate $f_\theta$;

$$\begin{aligned} f'_\theta(x) &= \frac{\theta}{\cosh^2(x/2) - \theta^2 \sinh^2(x/2)} \\ &= \frac{\theta}{1 + (1 - \theta^2) \sinh^2(x/2)} \qquad \forall x \in \mathbf{R}. \end{aligned} \tag{4.6}$$

The denominator in (4.6) is positive for all $x \in \mathbf{R}$ and increasing in $x$ for $x > 0$, so $f_\theta(x)$ is an increasing function of $x \in \mathbf{R}$ and a concave function for $x > 0$. Another consequence of (4.6) is that $f'_\theta(x)$ is an even function of $x$, whence $f_\theta(x)$, which vanishes at $x = 0$, is an odd function of $x$.

The denominator in (4.6) has the expansion $1 + (1 - \theta^2) x^2/4 + O(x^4)$ near 0, where the $O(x^4)$ term depends on $\theta$, but is a uniformly bounded multiple of $x^4$ for $\theta \in I$. Inverting and integrating, we see that the Taylor expansion of $f_\theta$ near 0 has the form,

$$f_\theta(x) = \theta x - \frac{\theta(1 - \theta^2)}{12} x^3 + O(x^5). \tag{4.7}$$

It remains to prove (4.5). Dividing that inequality by $\theta x$, inverting and squaring, shows that (4.5) is equivalent to $1+\kappa_1 x^2 \le (\frac{\theta x}{f_\theta(x)})^2 \le 1+\kappa_2 x^2$. In other words, we must verify that

$$\psi_\theta(x) := x^{-2}\left[\left(\frac{\theta x}{f_\theta(x)}\right)^2 - 1\right] \qquad \text{satisfies } \kappa_1 \le \psi_\theta(x) \le \kappa_2, \tag{4.8}$$

for all $x>0$ and $\theta \in I$, with some $\kappa_2 \ge \kappa_1 > 0$ that depend only on $I$. By (4.6), for $x>0$ we have $f_\theta(x) < \theta x$ so $\psi_\theta(x) > 0$. Therefore $\psi_\theta$ is uniformly bounded above and below by positive constants if $x$ and $\theta$ are restricted to compact intervals in $(0,\infty)$. Since

$$f_\theta(x) \to \log\left[\frac{1+\theta}{1-\theta}\right] \qquad \text{as } x\to\infty \text{ uniformly in } \theta\in I,$$

we deduce that $\psi_\theta(x)$ converges to a positive limit as $x\to\infty$, uniformly in $\theta \in I$. The expansion (4.7) implies that $\psi_\theta(x) \to \theta(1-\theta^2)/6$ as $x \to 0$. These considerations prove (4.8) and the lemma. □

Proof of Theorem 2.2. Specialize to plus boundary conditions. Thus we write $\mathbf{P}_v^{(N,+)}$ for $\mathbf{P}_v^{(N,\xi)}$ where $\xi \equiv +1$. Let

$$x_v^{(N)} = x_v^{(N,+)} = \log\left(\frac{\mathbf{P}_v^{(N,+)}(\eta(v)=+1)}{\mathbf{P}_v^{(N,+)}(\eta(v)=-1)}\right)$$

be the log-likelihood ratio of plus-to-minus at the root of the subtree $T^{(N)}(v)$. Note that with plus boundary conditions, all the $x_v^{(N)}$ are positive. Lemma 4.1 shows that

$$x_v^{(N)} = \sum_{v\to w} f_w(x_w^{(N)})$$

with $f_w = f_{\theta_w}$, as in (4.1).

Recall that the interaction strengths $J_v$ are in a bounded interval $[J_{\min}, J_{\max}] \subset (0,\infty)$, and $\beta$ is fixed. Therefore all the biases $\theta_v$ are in some bounded interval $I \subset (0,\infty)$. It follows from Theorem 3.2 and the inequalities in (4.5) that $x_o^{(N)}$ is bounded between $\mathrm{cap}_3(T^{(N)})/\kappa_2$ and $\mathrm{cap}_3(T^{(N)})/\kappa_1$ for all $N$. Taking decreasing limits finishes the proof of the theorem. □

We conclude this section with a discussion of the boundedness condition (1.1). Given any tree $T$ with associated interactions $\{J(e) : e \in E(T)\}$, a new tree $T'$ may be constructed by subdividing edges of $T$ according to the following scheme. Fix an $\varepsilon > 0$. Replace each edge $e$ with $\theta_e < \varepsilon$ by a series of $n$ edges $e(1), \ldots, e(n)$, with $\theta_{e(j)} = \theta_e^{1/n}$, where $n$ is the least integer making $\theta_e^{1/n}$ greater than $\varepsilon$.

From the error propagation description of the Ising measure, we see that the measure on $\{\pm 1\}^{V(T)}$, gotten by restricting the Ising measure on $T'$ to the vertices of $T$, coincides with the Ising measure on $T$. Distances in $T'$ no longer coincide with distances in $T$, but it is easy to see that the various definitions of phase transition in this article are unchanged if limits on $T'$ are taken with respect to distances in $T$. The associated resistor network to $T'$ may be described as follows. Each edge not subdivided retains the same resistance. A subdivided edge with resistance $R(e) = A/\theta_e$ is replaced by $n$ edges in series of resistances $A\theta_e^{-j/n}$ for $j = 1, \ldots, n$. Since $\theta_e^{1/n} < \varepsilon^{1/2}$, the effective resistance of these $n$ new edges in series is less than $1/(1 - \sqrt{\varepsilon})$ times the greatest resistance among them which is $A/\theta_e$. Thus the resistance of the new network is equal to the old resistance up to a bounded factor, and hence has capacity within a bounded factor of the original capacity. We conclude that no generality is lost by assuming $J_v$ to be bounded away from zero.

There is some generality lost in assuming $J_v$ to be bounded above, but for good reason, as shown by the following example. Let $T$ be a spherically symmetric tree with $|T_n| \approx n^\alpha 2^n$ for some $\alpha > 1/2$. As seen in Corollary 2.4, there is a phase transition on $T$ with constant interaction strength satisfying $\theta = 1/2$. Now replace each edge in generation $n$ by $n$ edges having $\theta_e = 2^{-1/n}$. The resistance of each new series of edges in generation $n$ is of order $n$ times the old resistance, so when $\alpha \leq 3/2$, the new tree has zero capacity. Thus the capacity criterion breaks down when the interaction strengths are allowed to have $\theta_v \to 1$, that is, $J_v \to \infty$.

**5. Free boundary conditions.** The question we ask in this section is: if you generate a configuration on $T^{(N)}$ from the free boundary measure, then look only at the boundary, do you have nonvanishing information about the root as $N \to \infty$? To formalize this, let $\xi$ be the random boundary configuration induced by the free measure $\mathbf{P}^{(N)}$ on configurations on all of $T^{(N)}$. In the notation of Lemma 4.1, let

$$X_v^{(N)} := x_v^{(N,\xi)}$$

be the log-likelihood ratio of plus-to-minus at $v$ given the boundary.

We want to know whether the $\mathbf{P}^{(N)}$ law of $X_o^{(N)}$ (the free law) converges weakly to the point mass at 0 as $N \to \infty$. Evans et al. (2000) showed that $X_o^{(N)}$ does not go to zero when $T$ has positive $L^2$ capacity with resistances given by (2.5). As mentioned in the Introduction, they, as well as Ioffe (1996b) have results in the other direction which leave the critical case open. We sharpen this by showing that zero capacity implies $X_o^{(N)} \xrightarrow{\mathcal{D}} 0$. The following statement is equivalent to Theorem 2.1.

THEOREM 2.1′. *Let $T$ be an infinite locally finite tree, rooted at $o$, with no leaves except possibly at $o$ and interaction strengths $J_v$ satisfying (1.1) and set $\theta_v = \tanh(\beta J_v))$. Suppose that $\mathrm{cap}_2(T) = 0$ with resistances as in (2.5). Then $X_o^{(N)}$ converges in law to 0.*

PROOF. By Lemma 4.1, when $|v| < N$,

$$X_v^{(N)} = \sum_{v \to w} f_w(X_w^{(N)}) \tag{5.1}$$

holds pointwise, with $f_w$ as in (4.1). To make use of this functional recursion, we will derive from it a system of real inequalities;

$$m_v^{(N)} \le \sum_{v \to w} \frac{\theta_w^2 m_w^{(N)}}{1 + \kappa m_w^{(N)}}. \tag{5.2}$$

The quantity $m_v^{(N)}$ will be an expectation of $X_v^{(N)}$ but it is not obvious what measure should be used to take the expectation. Define the measures $Q_v^{N+}$ (respectively, $Q_v^{N-}$) on the $\sigma$-field $\mathcal{F}_v^{(N)}$ of boundary values by letting

$$Q_v^{N+}(\xi) := \mathbf{P}_v^{(N)}(\eta : \eta_{|\partial T_v^{(N)}} = \xi | \eta(v) = +1)$$

be the conditional distribution of the free boundary given a plus at $v$ (respectively, given a minus at $v$). Define

$$m_v^{(N)} := \int X_v^{(N)} \, dQ_v^{N+} = -\int X_v^{(N)} \, dQ_v^{N-}.$$

The properties of the measures $Q_v^{N\pm}$ summarized in the following lemmas make these appropriate for the study of the free boundary.

LEMMA 5.1. *For any $v$ with $|v| < N$,*

$$Q_v^{N+} = \prod_{v \to w} \left[ \frac{(1+\theta_w)}{2} Q_w^{N+} + \frac{(1-\theta_w)}{2} Q_w^{N-} \right].$$

*In particular, the projection of $Q_v^{N+}$ onto boundary configurations on $T^{(N)}(w)$ is*

$$\frac{(1+\theta_w)}{2} Q_w^{N+} + \frac{(1-\theta_w)}{2} Q_w^{N-}.$$

LEMMA 5.2. *For any odd function $\phi$,*

$$\int \phi(X_v^{(N)}) \, dQ_v^{N+} = \int \phi(|X_v^{(N)}|) \tanh(|X_v^{(N)}|/2) \, d\mathbf{P}_v^{(N)}.$$

LEMMA 5.3. *There is a positive, continuous function $\kappa$ such that when $f_\theta$ is defined as in (4.1) with $\theta = \theta_v$, then*

$$\int f_\theta(X_v^{(N)})\, dQ_v^{N+} \le \theta \frac{\int X_v^{(N)}\, dQ_v^{N+}}{1 + \kappa(\theta) \int X_v^{(N)}\, dQ_v^{N+}}.$$

To finish the proof from these lemmas, use (5.1) and Lemma 5.1 to evaluate

$$\begin{aligned} m_v^{(N)} &= \sum_{v \to w} \int f_w(X_w^{(N)})\, dQ_v^{N+} \\ &= \frac{1}{2} \sum_{v \to w} \int f_w(X_w^{(N)})\, d((1+\theta_w) Q_w^{N+} + (1-\theta_w) Q_w^{N-}) \qquad (5.3) \\ &= \sum_{v \to w} \int \theta_w f_w(X_w^{(N)})\, dQ_w^{N+}. \end{aligned}$$

Apply Lemma 5.3 to see that this is, at most,

$$\sum_{v \to w} \frac{\theta_w^2 m_v^{(N)}}{1 + \kappa(\theta_v) m_v^{(N)}}.$$

By continuity of $\kappa(\theta)$ and the boundedness assumption (1.1), we arrive at (5.2). Theorem 3.2 now applies to show that $m_o^{(N)} \le \frac{\text{cap}_2(T^{(N)})}{\kappa}$ with resistances as in the hypothesis of the theorem. Hence $\text{cap}_2(T) = 0$ implies $m_o^{(N)} \to 0$ as $N \to \infty$. Finally, by Lemma 5.2 with $\phi(x) = x$, this implies that $X_o^{(N)} \xrightarrow{\mathcal{D}} 0$ as $N \to \infty$, finishing the proof. □

It remains to prove the lemmas. Lemma 5.1 is immediate from the Markov property.

PROOF OF LEMMA 5.2. We first compare $Q_v^{N+}$ to the boundary measure induced by the free measure $\mathbf{P}_v^{(N)}$. We claim that

$$\frac{dQ_v^{N+}}{d\mathbf{P}_v^{(N)}} = 1 + \tanh(X_v^{(N)}/2). \qquad (5.4)$$

Indeed, from Bayes' rule, one gets

$$\frac{dQ_v^{N+}}{d\mathbf{P}_v^{(N)}} = \frac{\mathbf{P}_v^{(N)}(\eta(v) = +1 | \mathcal{F}_v^{(N)})}{\mathbf{P}_v^{(N)}(\eta(v) = +1)}.$$

The denominator is $1/2$ by symmetry, while the numerator is $\exp(X_v^{(N)})/(1 + \exp(X_v^{(N)})) = (1 + \tanh(X_v^{(N)}/2))/2$ by the definition of $X_v^{(N)}$. This proves the claim. Now if $\phi$ is any odd function, then $\phi(x) = (\phi(x) - \phi(-x))/2$, and thus

$$\begin{aligned}\int \phi(X_v^{(N)})\, dQ_v^{N+} &= \int \frac{1}{2}(\phi(X_v^{(N)}) - \phi(-X_v^{(N)}))\, dQ_v^{N+} \\ &= \int (\phi(X_v^{(N)}) - \phi(-X_v^{(N)})) \frac{e^{X_v^{(N)}/2}}{e^{X_v^{(N)}/2} + e^{-X_v^{(N)}/2}}\, d\mathbf{P}_v^{(N)} \\ &= \int \phi(X_v^{(N)}) \frac{e^{X_v^{(N)}/2} - e^{-X_v^{(N)}/2}}{e^{X_v^{(N)}/2} + e^{-X_v^{(N)}/2}}\, d\mathbf{P}_v^{(N)}.\end{aligned}$$

The integrand is a product of two odd functions, whence it is an even function of $X_v^{(N)}$. Inserting absolute values yields the desired conclusion. □

Proof of Lemma 5.3. Abbreviate the notation by writing $X$ for $X_v^{(N)}$, $\mathbf{E}$ for integration against $\mathbf{P}_v^{(N)}$ and $\mathbf{E}_+$ for integration against $Q_v^{N+}$. First, for any $c > 0$, the product,

$$\mathbf{E}_+ f_\theta(X)(1 + c\mathbf{E}_+ X) = \mathbf{E}_+ f_\theta(X) + c(\mathbf{E}_+ f_\theta(X))(\mathbf{E}_+ X),$$

is equal, by Lemma 5.2, to the sum

$$\mathbf{E}[f_\theta(|X|)\tanh|X/2|] + \mathbf{E}[f_\theta(|X|)\tanh|X/2|)] \cdot \mathbf{E}[c|X|\tanh|X/2|].$$

Since the functions $f_\theta(x)\tanh(x/2)$ and $cx\tanh(x/2)$ are both nondecreasing on $[0, \infty)$, they are positively correlated functions of $|X|$ (under $\mathbf{P}_v^{(N)}$ or any other law), and hence

$$\begin{aligned}(\mathbf{E}_+ f_\theta(X))(1 + c\mathbf{E}_+ X) &\le \mathbf{E}[f_\theta(|X|)\tanh|X/2|] + \mathbf{E}[c|X| f_\theta(|X|)\tanh^2|X/2|] \\ &= \mathbf{E}[f_\theta(|X|)\tanh|X/2|(1 + c|X|\tanh|X/2|)].\end{aligned}$$

Recall that $\tanh(x) = x - \Theta(x^3)$. Refer to the Taylor expansion for $f_\theta = f_v$ in (4.7) to see that for $\kappa(\theta)$ sufficiently small, there is a range $x \in [0, \delta]$ for which

$$f_\theta(x)(1 + \kappa(\theta)x\tanh(x/2)) < \theta x. \tag{5.5}$$

Since $f_\theta$ is itself bounded and less than $\theta x - \varepsilon(\theta)x$ on $[\delta, \infty)$, we may choose $\kappa(\theta)$ smaller, if necessary, so that (5.5) holds for all $x \ge 0$. Clearly the choice of $\kappa$ can be made continuously in $\theta$. It follows that

$$(\mathbf{E}_+ f_\theta(X))(1 + \kappa(\theta)\mathbf{E}_+ X) \le \mathbf{E}[\theta|X|\tanh|X/2|] = \theta\mathbf{E}_+ X,$$

by Lemma 5.2. Dividing by $(1 + \kappa(\theta)\mathbf{E}_+ X)$ proves the lemma. □

**6. Spin-glasses.** Let $\mathbf{P}_v^{(N,\mathrm{sg})}$ denote the spin-glass measure $\mathbf{P}^{\mathrm{sg}}$ on configurations on the tree $T^{(N)}(v)$ (see Section 1 for definitions). Our object in this section is to determine when the conditional probability $\mathbf{P}_o^{(N,\mathrm{sg})}(\eta(o) = +1|\mathcal{F}^{(N)})$ converges in distribution to a point mass at 1/2 where $\mathcal{F}^{(N)} = \mathcal{F}_o^{(N)}$ is the $\sigma$-field generated by boundary values on $T^{(N)}$. By the Markov random field property (or by the definitions of $\mathbf{P}$ and $\mathbf{P}^{\mathrm{sg}}$), the measures $\mathbf{P}^{(N)}$ and $\mathbf{P}^{(N,\mathrm{sg})}$ agree when conditioned on the boundary, so the functions $X_v^{(N)}$ of the previous section compute conditional probabilities with respect to $\mathbf{P}^{(N,\mathrm{sg})}$. Thus our task is to see when $X_o^{(N)} \xrightarrow{\mathcal{D}} 0$ under the laws $\mathbf{P}^{(N,\mathrm{sg})}$.

THEOREM 2.5. *Let $T$ be an infinite, locally finite tree, rooted at $o$, with no leaves except possibly at $o$ and interaction strengths $J_v$ satisfying (1.1) and set $\theta_v = \tanh(\beta J_v)$. Then $X_o^{(N)} \xrightarrow{\mathcal{D}} 0$ under the spin-glass measure if and only if $\mathrm{cap}_2(T) = 0$ with resistances $R_v = \prod_{y \le v} \theta_y^{-2}$ as assigned in (2.5).*

PROOF. The structure of the proof is similar to that of Theorem 2.1. We begin with (5.1),

$$X_v^{(N)} = \sum_{v \to w} f_w(X_w^{(N)}).$$

Let $U_v^{(N)} := (X_v^{(N)})^2$ and

$$u_v^{(N)} := \int U_v^{(N)} \, d\mathbf{P}_v^{(N,\mathrm{sg})},$$

where the integrating measure in this case is just i.i.d. fair coin-flips on the boundary of $T^{(N)}(v)$. In place of Lemma 5.1 we have the observation that the random variables $X_w^{(N)}$ have mean zero and are independent as $w$ ranges over the children of a fixed $v$. Lemmas 5.2 and 5.3 are replaced by the following two lemmas. Define

$$g_v(x) := (f_v(\sqrt{x}))^2.$$

LEMMA 6.1. *For all $v$ and all $N > |v|$,*

$$\mathbf{E}(U_v^{(N)})^2 \le 3(\mathbf{E}U_v^{(N)})^2.$$

LEMMA 6.2. *There are continuous functions $\kappa_2(c,\theta_v) \ge \kappa_1(c,\theta_v) > 0$ such that for any random variable $V$ satisfying $\mathbf{E}V^2 \le c(\mathbf{E}V)^2$, one has*

$$h_2(\mathbf{E}V) \le \mathbf{E}g_v(V) \le h_1(\mathbf{E}V) \tag{6.1}$$

*with $h_i(x) = \theta_v^2 x/(1 + \kappa_i(c,\theta_v)x)$.*

From these two lemmas the proof is finished as follows. Let $\mathbf{E}$ denote expectation with respect to i.i.d. unbiased (spin-glass) boundary conditions. Since each $f_v$ is an odd function, the quantities $f(X_w^{(N)})$ are independent mean-zero as $w$ varies over the children of $v$, which gives rise to the recursive formula

$$\begin{aligned} u_v^{(N)} &= \mathbf{E}(X_v^{(N)})^2 \\ &= \mathbf{E}\Bigg(\sum_{v\to w} f_v(X_w^{(N)})\Bigg)^2 \\ &= \sum_{v\to w} \mathbf{E} f_v(X_w^{(N)})^2 \\ &= \sum_{v\to w} \mathbf{E} g_v(U_w^{(N)}). \end{aligned}$$

Apply Lemma 6.2 with $V=U_v^{(N)}$ and $c=3$ (obtaining the hypothesis from Lemma 6.1), to get

$$\sum_{v\to w} h_2(u_w^{(N)}) \le u_v^{(N)} \le \sum_{v\to w} h_1(u_w^{(N)}).$$

By continuity and the boundedness assumption (1.1), we may take $\kappa_i$ in the definition of $h_i$ to be constants independent of $v$. By Theorem 3.2 we see that $\lim_{N\to\infty} u_o^{(N)}$ is estimated up to a constant factor by $\mathrm{cap}_2(T)$ with resistances as stated in the hypothesis of the theorem. Since $X_o^{(N)}$ has mean zero and is bounded by $\sum_{o\to v}\log[(1+\theta_v)/(1-\theta_v)]$, it follows that the random variables $X_o^{(N)}$ converge in distribution to 0 if and only if their variances $u_o^{(N)}$ go to zero. This completes the proof of Theorem 2.5. □

It remains to prove Lemmas 6.1 and 6.2. Before proving Lemma 6.1, we record some preliminary facts.

LEMMA 6.3. *Suppose $f$ is a differentiable, weakly increasing and concave function on $[0,\infty)$ with $f(0)=0$. Then $x^2\circ f\circ\sqrt{x}$ is concave.*

PROOF. Let $\varphi(x)=f(x_0)+(x-x_0)f'(x_0)$ be the tangent line for $f$ at $x_0$. Concavity implies that $\varphi(x)\ge f(x)$ for all $x\ge 0$ and that $\varphi'(x_0)\le f(x_0)/x_0$. Thus $\varphi(x)=ax+b$ with $b\ge 0$, whence $x^2\circ\varphi\circ\sqrt{x}$ is a concave support function, lying above $x^2\circ f\circ\sqrt{x}$ with equality at $x_0^2$. We conclude that $x^2\circ f\circ\sqrt{x}$ is the minimum of a family of concave functions. □

LEMMA 6.4. *Let $g:[0,\infty)\to[0,\infty)$ be concave with $g(0)=0$, and let $Y$ be a nonnegative random variable with positive finite variance. Then*

$$\frac{\mathbf{E}[g(Y)^2]}{[\mathbf{E}g(Y)]^2}\le\frac{\mathbf{E}Y^2}{(\mathbf{E}Y)^2}. \tag{6.2}$$

PROOF. Let $Z=Y/\mathbf{E}Y$ and $h(z)=g(z\mathbf{E}Y)/\mathbf{E}(g(Y))$. Then $\mathbf{E}Z=\mathbf{E}h(Z)=1$, so there must exist $z_1,z_2>0$ such that $h(z_1)\ge z_1$ and $h(z_2)\le z_2$. Note that (6.2) is equivalent to $\mathbf{E}[h(Z)^2]\le\mathbf{E}[Z^2]$. We may assume that $h(z)$ is not identically equal to $z$, and thus, by concavity, there is a unique fixed point $x>0$ for which $h(x)=x$. For any $z\ge0$,

$$|h(z)-x|\le|z-x|,$$

and therefore,

$$\mathbf{E}[h(Z)^2]=\mathbf{E}(h(Z)-x)^2+2x-x^2\le\mathbf{E}(Z-x)^2+2x-x^2=\mathbf{E}Z^2,$$

proving the lemma. □

LEMMA 6.5. *For any nonnegative random variable $X\in L^4$, and any concave function $f$ with $f(0)=0$,*

$$\frac{\mathbf{E}f^4(X)}{(\mathbf{E}f^2(X))^2}\le\frac{\mathbf{E}X^4}{(\mathbf{E}X^2)^2}.$$

PROOF. by Lemma 6.3, the function $g:=x^2\circ f\circ\sqrt{x}$ is concave. Applying Lemma 6.4 to the function $g$ and the random variable $Y=X^2\in L^2$ gives

$$\frac{\mathbf{E}f^4(X)}{(\mathbf{E}f^2(X))^2}=\frac{\mathbf{E}g^2(Y)}{[\mathbf{E}g(Y)]^2}\le\frac{\mathbf{E}Y^2}{(\mathbf{E}Y)^2}=\frac{\mathbf{E}X^4}{(\mathbf{E}X^2)^2},$$

proving the lemma. □

REMARK. As noted by the referee, Lemmas 6.2, 6.4 and 6.5 are valid for quite general random variables; it would be interesting to apply them to more general situations.

PROOF OF LEMMA 6.1. Recall the definitions of $U_v^{(N)}$ and $u_v^{(N)}$ and define the fourth moment $s_v^{(N)}$:

$$U_v^{(N)}=(X_v^{(N)})^2;$$
$$u_v^{(N)}=\mathbf{E}U_v^{(N)};$$
$$s_v^{(N)}=\mathbf{E}(U_v^{(N)})^2=\mathbf{E}(X_v^{(N)})^4.$$

For any $v$, the random variables $\{f_w(X_w^{(N)}) : v \to w\}$ are independent with mean zero, so any monomial of these will have mean zero unless all exponents are even. The basic recursion (5.1) yields

$$\begin{aligned} u_v^{(N)} &= \mathbf{E}\biggl(\sum_{v\to w} f_w(X_w^{(N)})\biggr)^2 \\ &= \sum_{v\to w} \mathbf{E} f_w(X_w^{(N)})^2. \end{aligned}$$

Hence

$$(u_v^{(N)})^2 = \sum_{v\to w} (\mathbf{E} f_w(X_w^{(N)})^2)^2 + \sum_{v\to\{w,w'\}} 2\mathbf{E} f_w(X_w^{(N)})^2 \mathbf{E} f_{w'}(X_{w'}^{(N)})^2. \tag{6.3}$$

The fourth power expands similarly:

$$\begin{aligned} s_v^{(N)} &= \mathbf{E}\biggl(\sum_{v\to w} f_w(X_w^{(N)})\biggr)^4 \\ &= \sum_{v\to w} \mathbf{E} f_w(X_w^{(N)})^4 + \sum_{v\to\{w,w'\}} 6\mathbf{E} f_w(X_w^{(N)})^2 \mathbf{E} f_{w'}(X_{w'}^{(N)})^2. \end{aligned} \tag{6.4}$$

It is required to show that $s_v^{(N)} \le 3(u_v^{(N)})^2$.

Proceed by induction on $N - |v|$. First suppose $N - |v| = 1$ and that $v$ has $d$ children. Then $X_v^{(N)}$ is the sum of $d$ independent mean-zero random variables, each equal to $\pm\log(p/(1-p))$. In this case, $s_v^{(N)}/(u_v^{(N)})^2 = 3 - 2/d < 3$. Now suppose $N - |v| > 1$. By induction, $s_w \le 3u_w^2$ for each child $w$ of $v$. Applying Lemma 6.5, we see that for each such $w$,

$$\mathbf{E} f_w(X_w^{(N)})^4 \le 3(\mathbf{E} f_w(X_w^{(N)})^2)^2.$$

Plugging this into (6.4) and comparing with (6.3) shows that $s_v^{(N)} \le 3(u_v^{(N)})^2$, completing the induction. □

Proof of Lemma 6.2. We observed in the proof of Lemma 6.5 that $g_v$ is concave. For the upper bound, first note that

$$g_v(x) \le h(x) := \frac{\theta_v^2 x}{1 + \kappa(\theta_v) x},$$

for some $\kappa(\theta)$ is bounded above and below by positive constants for $\theta$ in a compact interval. The proof of this is the same as the proof of (4.5), using the Taylor expansion [that follows from (4.7)]

$$g_v(x) = \theta_v^2 x - \theta_v^2(1-\theta_v^2)x^2/6 + O(x^3),$$

together with boundedness and concavity of $g_v$. Jensen's inequality gives

$$\mathbf{E}g_v(V) \le \mathbf{E}h(V) \le h(\mathbf{E}V),$$

which proves the upper bound with $\kappa_1 = \kappa$.

For the lower bound, since $g_v(x) = \theta_v^2 x - O(x^2)$ near 0, we have $g_v(x) \ge \theta_v^2 x - \lambda x^2$ for some $\lambda$ and all $x$ in some interval $[0, \delta]$. Choosing $\lambda$ larger if necessary, we can ensure that $g_v(x) \ge \theta_v^2 x - \lambda x^2$ for all $x \ge 0$. Hence

$$\mathbf{E}g_v(V) \ge \theta_v^2 \mathbf{E}V - c\lambda(\mathbf{E}V)^2.$$

Choose $\delta(\theta_v) > 0$ so that the right-hand side is positive for $x \in (0, \delta(\theta_v))$. Choose $\kappa_2(\theta_v)$ so that

$$\frac{\theta_v^2 x}{1 + \kappa_2(\theta_v)x} \le [\theta_v^2 x - \lambda x^2] \wedge \frac{g_v(\delta/2)}{4c}.$$

This satisfies (6.1) when $\mathbf{E}V \le \delta$. But when $\mathbf{E}V > \delta$, then the hypothesis on $V$ implies that $\mathbf{P}(V > \delta/2) \ge 1/(4c)$ and therefore that $\mathbf{E}g_v(V) \ge g_v(\delta/2)/(4c)$. Hence (6.1) is valid for all $x \ge 0$. Together with the evident continuous dependence of $\kappa_i$ on $\theta_v$, this proves the lemma. □

**7. Concluding remarks.** Although we have in general no explicit probabilistic interpretation of $L^p$ capacities, in the case of integer values of $p$ there is a more probabilistic formulation. Positive $L^p$ capacity is equivalent to the existence of a probability measure $\mu$ on $\partial T$ such that $p$ independent paths picked from $\mu$ will coincide along a path of finite average resistance. This corresponds to the representation of $L^p$-energy as a $p$-fold integral over $\partial T$.

Finally, we note that other statistical mechanical models lead to recursions similar to (5.1) but with functions $f_v$ that are not necessarily concave. The Potts model with $1 < q < 2$ is essentially similar to the Ising model, but when $q > 2$, the functions $f_v$ are not concave and qualitatively different behavior arises. See Häggström (1996) for a discussion of this as pertains to the random cluster model.

Remark. Since the first draft of this paper was circulated in 1996, there have been many developments on the reconstruction problem, some of them influenced by that draft. As suggested by the referee, we summarize some of these developments here. Pemantle and Steif (1999) analyzed the Heisenberg model and other continuous-state models on general trees. They also introduced the important notion of "*Robust reconstruction*" where the boundary data is noisy. This notion was analyzed later in great generality by Janson and Mossel (2004). Census reconstruction on regular trees (where only the number of particles of each type on the boundary is given) was considered by Mossel and Peres (2003). A comprehensive survey of the area up to 2004 was written by Mossel (2004). A connection between reconstruction on trees and Glauber dynamics was found by Kenyon, Mossel and Peres (2001) [see also Berger et al. (2005)], and this theme was developed further by Martinelli, Sinclair and Weitz (2004). Notable progress on the reconstruction problem for the asymmetric Ising model was made by Borgs et al. (2006) and for the Potts model by Sly (2009). The arguments in Section 5 were extended to other boundary conditions in Ding, Lubetzky and Peres (2009) and used there to bound the relaxation time for Glauber dynamics at the critical temperature.

**Acknowledgments.** Much of the research presented here was performed at the Mittag Leffler Institute. We are grateful to E. B. Dynkin for telling us about $L^p$ capacities. We thank Manjunath Krishnapur, Gabor Pete, Antar Bandyopadhyay and the referee for helpful comments and corrections.

Department of Mathematics
University of Pennsylvania
209 South 33rd Street
Philadelphia, Pennsylvania 19104
USA
E-mail: pemantle@math.upenn.edu

Microsoft Research
1 Microsoft Way
Redmond, Washington 98052
USA
E-mail: peres@microsoft.com